\documentclass[final]{siamltex}
\usepackage{amsmath,amssymb}
\usepackage{cite}

\begin{document}

\title{Differential Calculus, Tensor Products and the Importance
of Notation}
\author{Jonathan H. Manton\thanks{
Control and Signal Processing Lab,
Department of Electrical and Electronic Engineering, 
The University of Melbourne, 
Victoria 3010 Australia.
Email: \texttt{j.manton@ieee.org}
}}

\maketitle

\pagestyle{myheadings}
\thispagestyle{plain}
\markboth{J. H. MANTON}{DIFFERENTIAL CALCULUS AND TENSOR PRODUCTS}

\begin{abstract}
An efficient coordinate-free notation is elucidated 
for differentiating matrix expressions and other functions between
higher-dimensional vector spaces.  This method of differentiation
is known, but not explained well, in the literature.  Teaching it
early in the curriculum would avoid the tedium of element-wise
differentiation and provide a better footing for understanding more
advanced applications of calculus.  Additionally, it is shown to
lead naturally to tensor products, a topic previously considered
too difficult to motivate quickly in elementary ways.
\end{abstract}

\newcommand{\reals}{\mathbb{R}}
\newcommand{\tr}[1]{\operatorname{tr}\left\{#1\right\}}

\section{Introduction}
\label{sec:intro}

The derivative of a function $f \colon \reals \rightarrow \reals$,
being a far-reaching concept, is taught early to students.
Higher-dimensional functions $f\colon \reals^n \rightarrow \reals^m$
can then be handled element-wise by computing the partial derivatives
of the components of $f = (f_1,\cdots,f_m)$. Yet an element-wise
approach to calculus is uninformative and tedious.  It does not
provide an inherently systematic way of differentiating matrix
expressions or functions between abstract vector spaces.

A known alternative exists.  It uses the notation $Df$ to denote
the (Fr\'echet) derivative~\cite[Chapter 8]{bk:Jost:analysis} of a
function $f\colon V \rightarrow W$ between normed vector spaces $V$ and
$W$.  The chain rule and product rule are typically written as
\begin{align}
\label{eq:chainrule}
D(f \circ g) &= (Df \circ g)\, Dg, \\
\label{eq:prodrule}
D(f\,g) &= (Df)\, g + f\, Dg.
\end{align}
Books remain silent on how to apply these rules. Even just
endeavouring to expand $D^2(f \circ g)$ leads to confusion though:
\begin{align}
D^2(f \circ g) 
    &= D\big( (Df \circ g)\, Dg \big) \\
    &= \big( D(Df \circ g) \big)\, Dg + (Df \circ g)\, D^2g  \\
    &= \big( (D^2f \circ g)\, Dg\big)\, Dg + (Df \circ g)\, D^2 g.
    \label{eq:faulty}
\end{align}
Taken literally, it makes no sense to multiply $D^2f \circ g$ with
$Dg$ twice. This confusion possibly caused an error in a well-regarded
book~\cite[p.\@ 3]{bk:Abraham:mappings} that was pointed out
in~\cite{MR0486377}.

This article propounds a minor modification of the $Df$ notation
that avoids such confusion, and exemplifies that the notation makes
differentiation easier, faster and more meaningful than working
exclusively with gradients, Jacobians and Hessians~\cite{bk:Magnus:matrix}.

The modification involves the tensor product $\otimes$.  Importantly,
$\otimes$ can be introduced merely as a formal symbol separating
the arguments of a function, and students can become familiar with
manipulating $\otimes$ as part of learning calculus.  Later, it can
be revealed that $\otimes$ is actually a tensor product that reduces
multi-linear maps to linear maps.  This pedagogic approach might
remove the difficulty students normally have with the concept of a
tensor product.

The details of how to use tensor products to simplify working with
derivatives are not readily found in the literature; no
mention is made in the following textbooks on differential
calculus~\cite[Chapter 2]{bk:Abraham:mfolds}, \cite[Chapter
8]{bk:Jost:analysis}, \cite[Chapter 4]{Zeidler:1986ij}, \cite[Chapter
5]{bk:Zorich:analysis_I}, nor in the following textbooks on
differential geometry~\cite[Chapter 1]{bk:Abraham:mappings},
\cite[Chapter I.2]{bk:Burke:diff}, \cite[Chapter I.3]{bk:Lang:diff_g},
\cite{Munkres:1991vu}.  Furthermore, Section~\ref{sec:rules}
illustrates that a formal treatment actually requires some care.

Interestingly, although calculus is often considered elementary, 
many aspects of it are
not elementary at all.  A plethora of articles exist on the
chain rule alone, including~\cite{dresden:deriv, MR0017784, MR0076845,
MR1531888, MR0486377, Rybakowski:frechet, Flanders:2001ka}.  The
existence of differentiable yet nowhere monotone
functions~\cite{Bruckner:1992hm}, while true, is far from obvious.
The history is not straightforward either; Fa\`a di Bruno was neither
the first to state nor prove the higher-order chain-rule formula
that bears his name~\cite{Johnson:bruno, MR2121322}. The present
article adds to this list by showing the traditional $Df$ notation
can trap the unwary.

\section{An Example in Matrix Space}
\label{sec:look}

The $Df$ notation provides a coordinate-free approach to differential
calculus.  It is first presented by example.

Consider $f(X) = \tr{X^T A X}$ where $\tr{}$ denotes trace, superscript
$T$ denotes transpose, and $A$ and $X$ are matrices of compatible
dimensions. It is known as the generalised Rayleigh
quotient~\cite{Helmke:1994ec} because the principal subspace of a
symmetric $A$ can be found by maximising $f(X)$ subject to the
normalising constraint $X^T X = I$.

Often the derivative of such a function $f$ is represented
by its Jacobian matrix whose $ij$-th element is the partial derivative
of $f$ with respect to the element $X_{ij}$ of $X$.  Evaluating
these partial derivatives from first principles is straightforward
but tedious: use $(AB)_{ij} = \sum_k A_{ik}B_{kj}$ twice and $\tr{Z}
= \sum_i Z_{ii}$ to obtain $f(X) = \sum_{ijk} X_{ji} A_{jk} X_{ki}$,
differentiate normally, and attempt to convert the answer back to
matrix form.

The following approach is considerably simpler.  Explanations follow in
subsequent sections.  Fix a matrix $Z$ of the same dimension as
$X$.  Then:
\begin{align}
f(X+tZ) - f(X) &= \tr{(X+tZ)^T A (X+tZ)} - \tr{X^T A X} \\
    &= \left( \tr{Z^T A X} + \tr{X^T A Z} \right) t
        + \left( \tr{Z^T A Z} \right) t^2. \label{eq:ftaylor}
\end{align}
Since derivatives represent linear approximations, (\ref{eq:ftaylor})
shows the derivative of $f$ at $X$ in the direction $Z$ is
$\tr{Z^T A X} + \tr{X^T A Z}$.  The meaning may not be clear yet,
but the calculation was simple!

The mapping $Z \mapsto \tr{Z^T A X} + \tr{X^T A Z}$
is linear:  if it sends $Z_1$ to $c_1$ and $Z_2$ to $c_2$
then it sends $\alpha Z_1 + \beta Z_2$ to $\alpha c_1 + \beta c_2$
for $\alpha,\beta \in \reals$.
This linear mapping is the (Fr\'echet) derivative of $f$.
\begin{align}
\label{eq:frechet}
Df(X)\cdot Z &= \tr{Z^T A X} + \tr{X^T A Z} \\
    &= \tr{Z^T (A + A^T) X}. \label{eq:ffre}
\end{align}
If required,
the Jacobian matrix can be read off as $(A+A^T)X$.

Treating $Z$ as a constant and differentiating (\ref{eq:ffre}) gives
\begin{equation}
\label{eq:D2f}
(D^2f(X) \cdot Z) \cdot T = \tr{Z^T(A+A^T)T}.
\end{equation}
The Hessian is $(A+A^T)$. The left-hand side of
(\ref{eq:D2f}) is more commonly written as $D^2f(X) \cdot (Z,T)$.

\section{First-Order Derivatives and Gradients}
\label{sec:fod}

The definition of the derivative
$f'(x) = \lim_{t \rightarrow 0} t^{-1} [ f(x+t) -
f(x) ]$ of a function $f\colon \reals \rightarrow
\reals$ extends in several ways to functions $f\colon U \rightarrow V$
between normed vector spaces $U$ and $V$.  The
reader may take, for concreteness, $U$ and $V$ to be scalars $\reals$,
vectors $\reals^n$ or matrices $\reals^{n \times m}$.  

One extension considers directional derivatives, reducing
to the case $g\colon \reals \rightarrow V$, $g(t) = f(x+tz)$ for
fixed $x,z \in U$, for which the same formula can be used:
\begin{equation}
\label{eq:gat}
D_z f(x) = \lim_{t \rightarrow 0} \frac{f(x+tz) - f(x)}{t}.
\end{equation}
If the limit exists for all $z$ then (\ref{eq:gat})
is called the G\^ateaux derivative of $f$ at $x$.

Another extension looks beyond (\ref{eq:gat}) and focuses on the
geometric meaning of $f'(x)$ as the gradient of the line of best
fit to the graph of $f$ at $x$.  This suggests defining the derivative
as the best linear approximation of $f$ at $x$. Precisely, fix $x$
and assume there exists a (continuous) linear function $A_x(z)$ such that
\begin{equation}
\label{eq:frech}
\lim_{z \rightarrow 0} \frac{\| f(x+z) - f(x) - A_x(z) \|}{\| z \|} = 0.
\end{equation}
Then $A_x$ is unique and is called the Fr\'echet derivative of $f$
at $x$, denoted $Df(x)$.  Sometimes, evaluation in a particular
direction is denoted using a dot, as in (\ref{eq:frechet}).  That
is, $Df(x) \cdot z = A_x(z)$.

The limit in (\ref{eq:frech}) must exist for \textit{all} sequences
$\{z_n\}_{n=1}^\infty$ with $z_n \rightarrow 0$.  Even if the mapping
$z \mapsto D_z f(x)$ in (\ref{eq:gat}) is linear for a fixed $x$,
the Fr\'echet derivative need not exist because 
(\ref{eq:frech}) may hold for sequences $z_n$ converging to the
origin along straight lines but not for sequences following certain
curved trajectories.  This occurs when the limit is not uniform
across straight lines: convergence to zero is fast along
some lines but arbitrarily slow along others.  (Appendix~\ref{sec:ce}
gives an example.)

Fr\'echet derivatives can be calculated by finding the directional
derivatives $D_zf(x)$ then verifying $A_x(z) = D_zf(x)$ satisfies
(\ref{eq:frech}).  Verification is unnecessary if the Fr\'echet
derivative is known to exist by other means.  The $f$ in Section~\ref{sec:look}
is a polynomial, hence its Fr\'echet derivative exists and can
be found using directional derivatives, either explicitly as in
(\ref{eq:ftaylor}) or, in more complicated situations, by using
truncated Taylor series approximations.  Of course, tables and rules
could be used instead.

If $f\colon U \rightarrow \reals$ is a scalar function then its gradient
at $x$ is defined with respect to an inner product.  This is often
forgotten because the Euclidean inner product is chosen without
mention in many textbooks.  In matrix space, the Euclidean inner
product is $\langle A,B \rangle = \tr{B^TA}$.  For a fixed matrix
$G$, $A(Z) = \langle G, Z \rangle$ is a linear functional, and every
linear functional can be written this way.  The gradient of $f$ at
$X$ is the matrix $G_X$ such that $Df(X) \cdot Z = \langle G_X, Z
\rangle$.

\section{Second-order Derivatives and Hessians}
\label{sec:sod}

The (Fr\'echet) derivative of a function $f\colon U \rightarrow V$
is $Df\colon U \rightarrow L(U;V)$ where $L(U;V)$ is the normed
vector space of (continuous) linear maps from $U$ to $V,$ with norm
the operator norm.
Applying $D$ to $Df$ yields the second-order derivative
$D^2f\colon U \rightarrow L(U;L(U;V))$.  A second-order derivative
requires not one, but two, directions: $(D^2f(X) \cdot T) \cdot Z$.
The right-hand side of (\ref{eq:Dforder}) interprets this as the
rate of change in the direction $T$ of the directional derivative
$Df(X) \cdot Z$.

To the letter of the law, $D^2f(X)$ is calculated from (\ref{eq:ffre})
as follows.  Working directly with $Df(X) \cdot Z$ is not allowed
because $Df(X)$ must be treated as an element of $L(U;V)$ when
computing $D^2f(X) \cdot T = D(Df)(X) \cdot T$.  By assuming the
Fr\'echet derivative exists, it suffices to work with directional
derivatives:
\begin{equation}
\label{eq:DDf}
D(Df)(X) \cdot T = \lim_{t \rightarrow 0} \frac{Df(X+tT) - Df(X)}{t}.
\end{equation}
For clarity, let $L_t = Df(X+tT) \in L(U;V)$.  For fixed $t$, both
$L_t - L_0$ and $(L_t - L_0)t^{-1}$ are linear operators in $L(U;V)$.
The vector space structure on $L(U;V)$ is that induced by pointwise
operations: $(L_t - L_0)t^{-1}$ evaluated at $Z$ is 
$(L_t \cdot Z - L_0 \cdot Z)t^{-1}$ by definition.  A sequence of linear
operators converges if and only if it converges pointwise (throughout,
all vector spaces are finite-dimensional for simplicity).  Thus, the right-hand
side of (\ref{eq:DDf}) can be determined pointwise:
\begin{align}
\left( \lim_{t \rightarrow 0} \frac{Df(X+tT) - Df(X)}{t} \right) \cdot Z
    &= \lim_{t \rightarrow 0} \frac{Df(X+tT) \cdot Z - Df(X) \cdot Z}{t} \\
    &= \tr{Z^T (A+A^T) T}. \label{eq:fsec}
\end{align}
In words, $D(Df)(X) \cdot T$ is the linear operator $Z \mapsto
\tr{Z^T (A+A^T) T}$.

A nominally different quantity is the derivative $Dg(X) \cdot T$
where $g(X) = Df(X) \cdot Z$ for a fixed $Z$.  Nevertheless, $Dg(X)
\cdot T = \tr{Z^T (A+A^T) T}$, the same as (\ref{eq:fsec}).  Indeed, the
pointwise vector space structure on $L(U;V)$ means
\begin{equation}
\label{eq:Dforder}
(D^2f(X)\cdot T) \cdot Z =
(D(Df)(X) \cdot T) \cdot Z = D(Df \cdot Z)(X) \cdot T.
\end{equation}
Therefore $D^2f$ can be calculated from $Df(X) \cdot Z$ by treating
$Z$ as a constant and differentiating with respect to $X$.
This is how (\ref{eq:D2f}) is obtained from (\ref{eq:ffre}).

The above notation is simple but cumbersome. Textbooks generally
drop the variables, writing the chain rule and product rule as
(\ref{eq:chainrule}) and (\ref{eq:prodrule}).  Without variables
though, deducing $(D^2(f \circ g)(X) \cdot T) \cdot Z$ from
(\ref{eq:faulty}) takes experience.

Including directions from the start reveals
\begin{align}
\label{eq:fcircg}
D(f \circ g) \cdot Z &= (Df \circ g) \cdot \big( Dg \cdot Z \big),   \\
\label{eq:fg}
D(f \cdot (g \cdot Z)) \cdot T &=
    (Df \cdot T) \cdot \big( g \cdot Z \big) +
    f \cdot \big( (Dg \cdot T) \cdot Z \big), \\
\begin{split}
\label{eq:DDfg}
(D^2(f \circ g) \cdot T) \cdot Z &=
    \big( (D^2f \circ g) \cdot (Dg \cdot Z) \big) \cdot (Dg \cdot Z) \\
        &\quad + (Df \circ g) \cdot \big( (D^2g \cdot T) \cdot Z \big).
\end{split}
\end{align}
Here, $X$ is omitted because it is simple enough to feed it in to
the terms requiring it. To be clear, $Df \circ g$ means evaluate
$Df$ at $g(X)$.

Neither approach is particularly friendly.  Omitting variables omits
important details while including variables is tedious; the reader
is invited to derive (\ref{eq:DDfg}) from either (\ref{eq:fcircg})
and (\ref{eq:fg}), or from (\ref{eq:chainrule}) and
\begin{equation}
\label{eq:Dfgcor}
D(f\,g) \cdot Z = (Df \cdot Z)\,g + f\, (Dg \cdot Z). 
\end{equation}

For scalar fields $f\colon U \rightarrow \reals$,
the unique linear operator $H_X$
satisfying $(D^2f(X) \cdot T) \cdot Z = \langle H_X \cdot T, Z
\rangle$ is the Hessian of $f$ at $X$.
The ordering is unimportant because $D^2f(X)$ is symmetric:
$(D^2f(X) \cdot T) \cdot Z = (D^2f(X) \cdot Z) \cdot T$ for all $Z$
and $T$.  When the Euclidean inner product is used, $H_X$ agrees
with what is called the Hessian matrix~\cite{bk:Magnus:matrix}.

\section{A Tensor Product Notation for Derivatives}
\label{sec:not}

Technically, the product $(Df)\,g$ in (\ref{eq:prodrule}) cannot
be formed; given $f\colon U \rightarrow L(V;W)$ and $g\colon U
\rightarrow L(Y;V)$, $Df$ maps into $L(U;L(V;W))$ whereas $g$ maps
into $L(Y;V)$.  The tensor product allows replacing $L(U;L(V;W))$
by $L(U \otimes V;W)$.  Equation (\ref{eq:prodrule}) should actually
be written
\begin{equation}
\label{eq:prod}
D(f\,g) = (Df)\,(I \otimes g) + f\,Dg
\end{equation}
where $I$ is the identity map.  This has the correctness of
(\ref{eq:Dfgcor}) and almost the same brevity as (\ref{eq:prodrule}).

The tensor product can be understood simply as directing variables to
their correct targets: the $g$ in $(Df)\,g$ blocks $Z$ from reaching
$Df$ when applied on the right, while the $I$ in $Df\,(I \otimes
g)$ allows the $Z$ through.  Although the direct sum could also
accomplish this, only the tensor product behaves
correctly under differentiation:
\begin{equation}
\label{eq:tpr}
D(f \otimes g) = (Df \otimes g) + (f \otimes Dg).
\end{equation}
In particular, (\ref{eq:prod}) can be differentiated again by using
$D(I \otimes g) = (DI \otimes g) + (I \otimes Dg) = (0 \otimes g)
+ (I \otimes Dg) = I \otimes Dg$.

If $f$ is itself a derivative then (\ref{eq:chainrule}) becomes
\begin{equation}
\label{eq:fgt}
D(f \circ g) = (Df \circ g)\,(Dg \otimes I).
\end{equation}
Following these rules gives
\begin{align}
D^2(f \circ g) &= (D^2f \circ g)\,(Dg \otimes I)\,(I \otimes Dg) +
    (Df \circ g)\,D^2g \\
\label{eq:D2ff}
    &= (D^2f \circ g)\,(Dg \otimes Dg) + (Df \circ g)\,D^2g, \\
\begin{split}
\label{eq:D3f}
D^3(f \circ g) &=
    (D^3f \circ g)\,(Dg \otimes Dg \otimes Dg) \\ &\quad + 
    (D^2f \circ g)\,\big[ (D^2g \otimes Dg) + 2(Dg \otimes D^2g) \big] +
    (Df \circ g)\,D^3g.
\end{split}
\end{align}
The remainder of this section gives the intermediate steps.
Section~\ref{sec:rules} presents a formal description of the notation.

Start with $D(f \circ g) = (Df \circ g)\,Dg$.  
Differentiate to
get $D(Df \circ g)\,(I \otimes Dg) + (Df \circ g)\,D^2g$.  This
time, (\ref{eq:fgt}) is required: $D(Df \circ g) = (D^2f \circ
g)\,(Dg \otimes I)$.  Tensor products of linear maps satisfy the
rule $(A \otimes B)\,(C \otimes D) = (AC \otimes BD)$.  Therefore,
$(Dg \otimes I)\,(I \otimes Dg) = Dg \otimes Dg$.

To obtain (\ref{eq:D3f}), first apply the product rule (\ref{eq:prod})
to the two additive terms in (\ref{eq:D2ff}).  Note $D(Dg \otimes
Dg) = (D^2g \otimes Dg) + (Dg \otimes D^2g)$.  At this point,
\begin{align}
\begin{split}
\label{eq:D3ftemp}
D^3(f \circ g) &= (D^3f \circ g) \, (Dg \otimes I)
        \, (I \otimes (Dg \otimes Dg)) \\ &\quad
    + (D^2f \circ g)\,\big[ (D^2g \otimes Dg) + (Dg \otimes D^2g) \big] \\ &\quad
    + (D^2f \circ g) \, (Dg \otimes I) \, (I \otimes D^2g)
    + (Df \circ g)\,D^3g.
\end{split}
\end{align}
The first $I$ in (\ref{eq:D3ftemp}) acts on $U \otimes U$ whereas
the second acts on $U$.  Regardless, it is agreeable to equate $(Dg
\otimes I)\,(I \otimes (Dg \otimes Dg))$ with $Dg \otimes (Dg \otimes
Dg) = Dg \otimes Dg \otimes Dg$, and (\ref{eq:D3f}) readily follows
from (\ref{eq:D3ftemp}).

\section{Formal Description}
\label{sec:rules}

The notation used in Section~\ref{sec:not} is derived below.  Some
intricacies appear but go unnoticed in practice.  The notation
simplifies the differentiation by hand of abstract expressions,
such as when seeking bounds like the one in (\ref{eq:bound}).  It
is generally not needed for differentiating specific functions; see
Section~\ref{sec:look}.

All spaces are finite-dimensional vector spaces.  Basic properties
of tensor products are used~\cite{bk:Yokonuma:tensor}.  The main
principle is that canonical isomorphisms of vector spaces can be
applied freely because they  essentially commute with the Fr\'echet
derivative.

\newcommand{\D}[3][]{\ensuremath{\bar D^{#1}_{#2}\!\left(#3\right)}}

Given $f\colon U \rightarrow V$, define $\bar D^k f\colon U \rightarrow
L(U \otimes \cdots \otimes U;V)$ by
\begin{equation}
\bar D^kf(X) \cdot (Z_1 \otimes \cdots \otimes Z_k)
    = ((D^kf(X) \cdot Z_1) \cdots) \cdot Z_k.
\end{equation}
For $g\colon U \rightarrow L(V;W)$, define $\D[k]{V}{g}\colon U \rightarrow
L(U \otimes \cdots \otimes U \otimes V;W)$ by
\begin{equation}
\D[k]{V}{g} \cdot (Z_1 \otimes \cdots \otimes Z_k \otimes T)
    = (((D^kg(X) \cdot Z_1) \cdots) \cdot Z_k) \cdot T.
\end{equation}
Although $\bar D^2f \neq \bar D(\bar Df)$, they agree up to
a canonical linear isomorphism.  In fact, $\bar D^kf = \D[k-1]{U}{\bar
Df}$.  For all intents and purposes, $\D[2]{V}{g}$ agrees with $\D{U
\otimes V}{\D{V}{g}}$ because applying the canonical identification
$U \otimes (U \otimes V) \cong U \otimes U \otimes V$ in practice
simply means omitting a pair of brackets.

Given $g\colon U \rightarrow L(V;W)$ and $h\colon U \rightarrow L(W;Y)$,
the product rule is
\begin{equation}
\label{eq:p1}
\D{V}{h \, g} = \D{W}{h} \, (I_U \otimes g) + h \, \D{V}{g}
\end{equation}
where $I_U\colon U \rightarrow U$ is the identity map and $I_U \otimes
g$ is a tensor field over $U$ whose value at $X \in U$ is $I_U
\otimes (g(X)) \in L(U \otimes V; U \otimes W)$.

Related is the application of a linear map to a vector; given $f\colon
U \rightarrow W$ then
\begin{equation}
\label{eq:Dhdotf}
\bar D(h \cdot f) = \D{W}{h} \, (I_U \boxtimes f) + h \, \bar Df
\end{equation}
where $(I_U \boxtimes f)(X)$ is the linear map $Z \mapsto (Z \otimes
f(X))$.  Later, by minor abuse of notation, $\otimes$ will replace
$\boxtimes$. Since $L(U;V) \boxtimes W = L(U;V \otimes W) \cong
L(U;V) \otimes W$, both $\boxtimes$ and $\otimes$ behave essentially the
same way when differentiated.

For $d\colon V \rightarrow L(W;Y)$, $e\colon U \rightarrow V$ and $f\colon V
\rightarrow W$, the two chain rules are
\begin{align}
\label{eq:c1}
\bar D(f \circ e) &= (\bar Df \circ e)\,\bar De, \\
\label{eq:c2}
\D{W}{d \circ e} &= (\D{W}{d} \circ e)\,(\bar De \otimes I_W)
\end{align}
where $I_W\colon W \rightarrow W$ is the identity map.

For $e\colon U \rightarrow V$ and $f\colon U \rightarrow W$, the tensor product rule is
\begin{equation}
\label{eq:t1}
\bar D(e \otimes f) = (\bar D e \boxtimes f) + (e \boxtimes \bar Df)
\end{equation}
where $\boxtimes$ combines a vector and a linear map to form a linear map,
as in (\ref{eq:Dhdotf}).
For $g\colon U \rightarrow L(V;W)$ and $h\colon U \rightarrow L(C;Y)$, the tensor
product rule is
\begin{equation}
\label{eq:t2}
\D{V \otimes C}{g \otimes h} = (\D{V}{g} \otimes h) + (g \otimes \D{C}{h}).
\end{equation}
For $e\colon U \rightarrow V$ and $h\colon U \rightarrow L(C;Y)$,
\begin{equation}
\label{eq:t3}
\D{C}{e \boxtimes h} = (\bar De \otimes h) + (e \boxtimes \D{C}{h}).
\end{equation}

If the codomains of all functions are spaces of linear maps then
the situation is particularly simple; (\ref{eq:p1}), (\ref{eq:c2})
and (\ref{eq:t2}) suffice.  This is the typical situation when
computing higher-order derivatives because the codomain of the
derivative of a function is a space of linear maps.  It is possible
to reduce to this situation by replacing $f\colon U \rightarrow
W$ with $\tilde f\colon U \rightarrow L(\reals;W)$ where $f(X) = \tilde
f(X) \cdot 1$.  The $\cdot 1$ can be removed, the derivatives
calculated, and the $\cdot 1$ applied at the very end.  This explains
the similarity of (\ref{eq:p1}) and (\ref{eq:Dhdotf}), and of
(\ref{eq:t1}), (\ref{eq:t2}) and (\ref{eq:t3}).

In practice, it is easier to replace $\boxtimes$ by $\otimes$ than
replace $f$ by $\tilde f$.  No confusion arises because $\boxtimes$
and $\otimes$ behave the same way with respect to addition,
multiplication and differentiation.

The subscripts on $\bar D$ used in (\ref{eq:p1})--(\ref{eq:t3})
merely keep all derivatives in a consistent form and can be dropped.
When computing higher-order derivatives recursively, to account for
$\bar D^k$ differing from $\bar D^{k-1}$ by a linear isomorphism,
it is only necessary to remove any remaining brackets in tensor
products at the end of each step, e.g., replace $Dg \otimes (Dg
\otimes Dg)$ by $Dg \otimes Dg \otimes Dg$ in (\ref{eq:D3ftemp}).

Once $\boxtimes$ is replaced by $\otimes$ and the subscripts dropped
on $\bar D$, the rules collapse to those in Section~\ref{sec:not}.

\section{Discussion}
\label{sec:ext}

Attention has been restricted to finite-dimensional vector spaces.
In principle, the results remain valid in infinite dimensions, but
a subtlety is that tensor products are not uniquely defined on
Banach spaces; different choices of norms, and hence completions
with respect to that norm, are possible~\cite{Ryan:2002wk}.

Bounds on the (operator) norms of derivatives are important in a
number of contexts, including for analysing the convergence of
iterative algorithms in numerical analysis. The formal treatment
in Section~\ref{sec:rules} justifies the calculation
\begin{align}
\| D^2(f \circ g) \|
&\leq \| D^2f \circ g \|\, \| Dg \otimes Dg \| + \| Df \circ g \|\, \| D^2g\| \\
\label{eq:bound}
&\leq \| D^2f \circ g \|\, \| Dg \|^2 + \| Df \circ g \|\, \| D^2g\|. 
\end{align}

An alternative derivation to that of Section~\ref{sec:rules} could
be based on treating $\otimes$ as a formal symbol used to direct
variables to their correct targets and developing a mechanical
calculus.  This would follow the course of building a class $\Omega$
of allowable expressions, explaining how $D$ is applied to members
of this class, and verifying the class is algorithmically closed
under $D$.

\section{Relevance}
\label{sec:rel}

This section discusses several situations where coordinate-free
differentiation simplifies matters.

The opening sentence of~\cite{MR1531888} asserts that formulae for
differentiating composite functions are simple only in the case of
first-order derivatives, where the chain rule applies.  Yet requiring
the first few derivatives of a composite function is a common
occurrence, such as requiring for a Newton method the first two
derivatives of a cost function $f\colon \reals^{n \times m} \rightarrow
\reals$ given as the composition $f = g \circ h$.

The derivatives of $f$ may be calculated from the formula
in~\cite{MR0486377} or~\cite{Rybakowski:frechet}, but unless one
is already familiar with the notation in (3.1)--(3.2)
of~\cite{Rybakowski:frechet}, the formula may be difficult to apply
quickly with confidence.  Perhaps less daunting is the formula
labelled the chain rule for Hessian matrices, appearing in~\cite[p.
110]{bk:Magnus:matrix}.  Yet this involves the Kronecker
product and may not yield a parsimonious description of the
second-order derivative.  Furthermore, it is inapplicable if the
domain of $f$ is not a matrix space.

The author's preferred choice is using (\ref{eq:chainrule}),
(\ref{eq:prod}) and (\ref{eq:fgt}) to differentiate $f = g \circ
h$ twice by hand.  No complicated formulae are involved, and the
same basic rules apply regardless of the actual domains of $g$ and
$h$.  Moreover, the standard rules for manipulating norms hold,
hence bounds such as (\ref{eq:bound}) are readily obtained.

The practical relevance of composite cost functions $f = g \circ
h$ on interesting domains includes the theory and practice of
optimisation on manifolds~\cite{Edelman:1998ei, Manton:2002hg,
Manton:2012vk}, where $g$ is a cost function on the manifold itself
and $h$ is used to ``pull'' the cost function back locally to a
function on the tangent space.  Matrix manifolds such as the Stiefel
manifold occur naturally in signal processing and the method
exemplified in Section~\ref{sec:look} is often the easiest way of
differentiating functions on matrix manifolds.

A statistician wishing to estimate from data the entries of a
symmetric matrix may be lead to studying cost functions whose domains
are symmetric matrices.  The coordinate-free framework handles this
effortlessly: $Df(X) \cdot Z$ is defined exactly as before, where
$X$ and $Z$ are symmetric matrices.

The function $f(X) = \log \det X$ of an invertible matrix $X$ is
encountered in various situations, including in relation to maximum
entropy methods~\cite{Pavon:2013kh}.  (The article~\cite{Pavon:2013kh}
itself implicitly advocates the coordinate-free approach to
differentiation because the coordinate-free approach makes transparent
the underlying geometry.) Differentiating $f(X) = \log \det X$
element-wise~\cite[pp.  149--151]{bk:Magnus:matrix} is tedious and
uninformative; as if by magic, a pleasing expression results.  By
comparison, a coordinate-free derivation is easily written down and
remembered.  The starting point is $\det(I+tZ) = I + \tr{Z} t +
\cdots$ which can be derived by writing $Z$ in Jordon canonical
form $Z = P^{-1} J P$ and expressing the determinant of an
upper-triangular matrix as the product of the diagonal elements.
Precisely,
\begin{align}
\det(I+tZ) &= \det\left(P^{-1} (I+tJ) P\right) \\
    &= \det(I+tJ) \\ &= \prod_i (1+tJ_{ii}) \\ &= 1 + t \sum_i
    J_{ii} + \cdots \\ &= 1 + t \tr{J} + \cdots .
\end{align}
Therefore,
\begin{align} f(X+tZ) &= \log
\det\left(X(I+tX^{-1}Z)\right) \\
    &= \log \det X + \log \det \left( I + t X^{-1}Z \right) \\ &=
    \log \det X + \log \left( 1 + t \tr{X^{-1} Z} + \cdots \right)
    \\ &= \log \det X + t \tr{X^{-1}Z} + \cdots
\end{align}
from which it follows immediately that $Df(X) \cdot Z = \tr{X^{-1}Z}$.

\section{Conclusion}

Derivatives of matrix expressions arise frequently in the applied
sciences~\cite{bk:Magnus:matrix}. The traditional element-wise
approach is tedius and uninformative compared with the $Df$ notation
(Section~\ref{sec:look}).  It is incongruous with the ease with
which $Df$ can be taught that it is not as widely used as it
profitably could be.

One could speculate the downfall of the $Df$ notation is the
difficulty encountered when repeatedly applying the chain and product
rules (Section~\ref{sec:intro}).  This difficulty is eliminated
by adopting the modified notation introduced in Section~\ref{sec:not}.

This modified notation is advocated to be taught to students early
in the curriculum. The tensor product $\otimes$ appearing in the
notation can be treated merely as a formal symbol separating arguments
to functions and which is differentiated analogously to the product
rule, hence the $\times$ in $\otimes$.  Furthermore, the notation
is pedagogically interesting as an elementary yet genuine application
of the tensor product.

\section*{Acknowledgement}

The author gratefully acknowledges conversations with Prof.\@
Louis Rossi that led to an improved presentation and the inclusion
of Section~\ref{sec:rel}.

\appendix
\section{A Counterexample}
\label{sec:ce}

Even if the directional derivatives (\ref{eq:gat}) exist and fit
together linearly, the derivative (\ref{eq:frech}) need not exist.
Simple counterexamples are known. It is nevertheless insightful to
derive a counterexample from first principles.

Consider the region in $\reals^2$ between the parametrised curves
$t \mapsto (t,0)$ and $t \mapsto (t,4t^2)$.  If $f\colon \reals^2
\rightarrow \reals$ is zero on and outside the boundary of the
parametrised region, it will be initially zero for a short distance
on every ray emanating from the origin, that is, $D_zf(0,0) = 0$.
If $f$ additionally satisfies $f(t,2t^2) = t$ --- a ridge of height
$t$ running along the curve $t \mapsto (t,2t^2)$ in the middle of
the aforementioned region --- then the derivative (\ref{eq:frech})
of $f$ at the origin would not exist.

Rigorously, let $\alpha\colon \reals \rightarrow [0,1] \subset
\reals$ be a smooth bump function that is zero outside the open
interval $(1,3) \subset \reals$ and which satisfies $\alpha(2)=1$.
Let $f(x,y) = x\,\alpha(yx^{-2})$ if $x \neq 0$ and $f(x,y) = 0$
otherwise.  Away from the $y$-axis, $f$ is smooth and \emph{a
fortiori} continuous.  If $(x_n,y_n) \rightarrow (0,y)$ then $|
f(x_n,y_n) | \leq | x_n | \rightarrow 0$, proving $f$ is everywhere
continuous. If $g(t) = f(at,bt)$ then there exists an $\epsilon >
0$ such that $g(t) = 0$ for $|t| < \epsilon$, that is, all directional
derivatives at the origin are zero. This means that if $Df(0,0)$
exists it must be zero, yet the sequence $(x_n,y_n) = (n^{-1},
2n^{-2})$ for $n=1,2,\cdots$ is such that $| f(x_n,y_n) | \, \|
(x_n,y_n) \|^{-1} \rightarrow 1 \neq 0$. (The Euclidean norm has been
used.)

\bibliographystyle{siam}
\bibliography{diffcalc}

\end{document}